\documentclass[12pt,a4paper]{article}
\usepackage{amssymb,amsfonts,amsthm,amsmath}
\usepackage{graphicx}
\usepackage{placeins}

\def\hang{\hangindent\parindent}
 \def\rf{\par\noindent\hang}
\begin{document}

\baselineskip=22pt

\begin{center} \large{{\bf THE PERFORMANCE OF A TWO-STAGE ANALYSIS OF ABAB/BABA CROSSOVER TRIALS}}
\end{center}

\bigskip

\noindent Running Title:
TWO-STAGE ANALYSIS OF CROSSOVER TRIALS

\bigskip

%\vskip .5in

\begin{center}
\large{{\bf PAUL KABAILA$^{1 *}$}} \ \normalsize{AND} \
\large{{\bf MATTHEW VICENDESE$^1$}}
\end{center}

\begin{center}
\large{\it La Trobe University}
\end{center}

%\bigskip
%\medskip

\vspace{12cm}

%\normalsize

\noindent $^*$ Author to whom correspondence should be addressed.

\noindent $^1$ Department of Mathematics and Statistics, La Trobe University, Victoria 3086, Australia.

\noindent e-mail: P.Kabaila@latrobe.edu.au

\noindent Facsimile: 3 9479 2466

\noindent Telephone: 3 9479 2594

\newpage

\begin{center}
\large{{\bf Summary}}
\end{center}

\noindent Freeman has considered the following two-stage procedure for finding
a confidence interval for the treatment difference theta, using data from an AB/BA
crossover trial. In the first stage, a preliminary test of the null hypothesis
that the differential carryover is zero, is carried out. If this hypothesis is accepted
then the confidence interval for theta is constructed assuming that the differential
carryover is zero. If, on the other hand, this hypothesis is rejected then this confidence
interval is constructed using only data from the first period.
Freeman has shown that this confidence interval has minimum coverage
probability far below nominal. He therefore concludes that this confidence interval should not
be used. In the present paper, we analyse the performance of a similar two-stage
procedure for an ABAB/BABA crossover trial.
This trial differs in very significant
ways from an AB/BA crossover trial, including the fact that for
an ABAB/BABA crossover trial there is an unbiased estimator of
the differential carryover that is unaffected by between-subject variation.
Despite these great differences, we arrive at the same conclusion as Freeman.
Namely, that the confidence interval resulting from the two-stage procedure
should not be used.

\bigskip


\newpage

\begin{center}
{\large{\bf 1. Introduction}}
\end{center}

 Consider a two-treatment two-period crossover trial, with
continuous responses.
%responses that are continuous random
%variables.
The purpose of this trial is to find a $1-\alpha$ confidence interval for the difference $\theta$
in the effects of two treatments, labelled A and B. Subjects are randomly allocated to either group 1
or group 2. Subjects in group 1 receive treatment A in the first period and then
%(after a washout
%period)
receive treatment B in the second period. Subjects in group 2 receive treatment B in the first
period and then
%(after a washout period)
receive treatment A in the second period. This trial is called
an AB/BA trial. To deal with the possibility of non-zero differential carryover, it was suggested (starting with
Grizzle, 1965, 1974 and endorsed by Hills \& Armitage, 1979 and
Armitage \& Hills, 1982) that the following two-stage procedure
be used. In the first stage, a preliminary test of the null hypothesis that the differential carryover is
zero (against the alternative that it is non-zero) is carried out. If this null hypothesis is accepted
then the confidence interval for $\theta$ is constructed to have nominal coverage
%probability
$1-\alpha$,
assuming that there is no differential carryover.
If, on the other hand, this null hypothesis is rejected then this confidence interval is constructed using
only data from the first period (since this is unaffected by carryover). As pointed out by Freeman (1989),
accepting this null hypothesis is not equivalent to concluding that the differential carryover is exactly zero.
Freeman (1989) shows that the confidence interval interval resulting from this two-stage procedure has minimum
coverage probability far below $1-\alpha$, demonstrating that this confidence interval should not be used.
%is completely inadequate.
%The conclusion from this is that, when analysing data from an AB/BA crossover trial, this two-stage
%procedure should not be used.
Senn (2006) states ``In my opinion the most important paper on cross-over trials in the 25 years of
{\it Statistics in Medicine} is Peter Freeman's paper''
% ... However, this was not only the most important
%paper on cross-over trials in {\it Statistics in Medicine}, but also the most important paper on this theme
%anywhere''.

What is the performance of this type of two-stage procedure for other crossover designs?
Jones \& Kenward (2003, pp. 123--125) analyse the performance of this type of procedure for Balaam's design.
This analysis makes the following two assumptions. The first assumption is that any carryover
from a treatment in a given period is only into the next period, and not beyond
(``first-order carryover'' model). The second assumption is that the carryover from
one period into the next period is determined only by the treatment applied in the first
period and not the treatment applied in the second period. Thus, for example, according
to this second assumption the carryover from treatment A into the next period is the same, irrespective
of whether the treatment in the next period is A or B.
This assumption has rightly been criticized
as being unrealistic by Fleiss (1986, 1989), Senn \& Lambrou (1998) and Senn (2001, 2002, 2005). This severely limits
the applicability of the analysis of Jones \& Kenward (2003) of this type of procedure for
Balaam's design.

In the present paper we consider an ABAB/BABA crossover trial.
Subjects are randomly allocated to either group 1
or group 2. Subjects in group 1 receive treatments A, B, A and B in the first, second, third and fourth
periods respectively. Subjects in group 2 receive treatments B, A, B and A in the first, second, third and fourth
periods respectively. We assume that any carryover
from a treatment in a given period is only into the next period, and not beyond.
However,
our analysis of this trial
does {\sl not} require us to assume that the carryover from
one period into the next period is determined only by the treatment applied in the first
period and not the treatment applied in the second period.
This is because we never need to consider the
carryover of a treatment from one period into the next period for which the same treatment is
applied. Two major differences between the AB/BA and ABAB/BABA trials are the following.
For an ABAB/BABA trial:
\begin{enumerate}

\item[(i)] There is an unbiased estimator (which is unaffected by differential carryover) of $\theta$
that has the following properties. It is unaffected by the between-subject variation. Also, it is
obtained without ignoring all of the data from periods 2, 3 and 4. This is the estimator $\hat \Theta$
described in Section 2.

\item[(ii)] There is an unbiased estimator of the differential carryover that is unaffected by the between-subject
variation. This is the estimator $\hat \Psi$ described in Section 2.

\end{enumerate}

There are two arguments against the adoption of $\hat \Theta$ as the standard estimator of $\theta$.
Firstly, as shown in Appendix A, this estimator is inefficient by comparison with the usual estimator
of $\theta$ based on data from a completely randomized design, using the same number of measurements of 
the response, unless a restrictive condition
holds. Secondly, there is an estimator of $\theta$, which we denote by $A$ and describe in Section 2,
that is much more efficient than $\hat \Theta$, when the differential carryover is zero. We view
$\hat \Theta$ as the analogue for an ABAB/BABA design of the estimator of $\theta$ constructed
using only data from the first period of an AB/BA design.

To deal with the possibility of non-zero differential
carryover, it is tempting to consider the use of the following two-stage procedure.
In the first stage a preliminary test of the null
hypothesis that the differential carryover is zero (against the alternative that it is non-zero)
is carried out. If this null hypothesis is accepted
then the confidence interval for $\theta$ is constructed using the estimator $A$ (described in Section 2)
and having nominal coverage
%probability
$1-\alpha$,
assuming that there is no differential carryover.
If, on the other hand, this null hypothesis is rejected then this confidence interval is constructed
to have nominal coverage $1-\alpha$, using
the estimator $\hat \Theta$ (described in Section 2) that is
based on data from all 4 periods.
This two-stage procedure is described in detail in Section 2.

A computationally-convenient formula for the coverage probability of the confidence interval
that results from this procedure is presented in Section 2.
In Section 3 we numerically evaluate the coverage properties of this confidence interval.
We show that this confidence interval has minimum
coverage probability far below $1-\alpha$, demonstrating that this confidence interval should
not be used. The coverage probability of this confidence interval depends only on the scaled
differential carryover. This is in sharp contrast to the coverage
probability of the confidence interval resulting from the two-stage procedure applied
to data from an AB/BA trial, found by Freeman (1989), which
depends on both the scaled differential carryover and the ratio
(error variance)/(subject variance).

Beginning with the work of Freeman (1989), the literature on the effect of preliminary model
selection (using, for example, hypothesis tests or minimizing a criterion such as AIC or
Mallows's $C_P$) on confidence intervals has grown steadily. This literature is reviewed by
Kabaila (2009). It is commonly the case that preliminary model selection has a highly detrimental effect on the
coverage probability of these confidence intervals.
However, each case needs to be considered individually on its merits.

\newpage

%\bigskip

\begin{center}
{\large\textbf{2. The two-stage analysis of ABAB/BABA trials under  consideration}}
\end{center}

We assume the following model for the ABAB/BABA trial. This model
is similar to the model for an AB/BA crossover trial
put forward by Grizzle (1965), as described by Grieve (1987).
Let $n_1$ and $n_2$ denote the number of subjects in group 1 and group 2 respectively.
Also let $Y_{ijk}$ be the response of the $j$th subject in the $i$th group and the $k$th period
($i=1,2$;  $j=1,\dots, n_i$; $k=1,2,3,4$). The model is
\begin{equation}
\label{model}
Y_{ijk} = \mu + \xi_{ij} + \pi_k + \phi_{\ell} + \lambda_q + \varepsilon_{ijk}
\end{equation}
where \newline
\phantom{1} \ \ \ $\mu$ is the overall population mean \newline
\phantom{1} \ \ \ $\xi_{ij}$ is the effect of the $j$th subject in the $i$th group \newline
\phantom{1} \ \ \ $\pi_k$ is the effect of the $k$th period \newline
\phantom{1} \ \ \ $\phi_{\ell}$ is the effect of the $\ell$th treatment \newline
\phantom{1} \ \ \ $\lambda_q$ is the residual effect of the $q$th treatment \newline
\phantom{1} \ \ \ $\varepsilon_{ijk}$ is the random error \newline
Note that both $\ell$ and $q$ are determined by the group $i$ and the period $k$.
This model is described in less abbreviated form in Appendix A.
We assume that the $\xi_{ij}$ and $\varepsilon_{ijk}$ are independent and that the
$\xi_{ij}$ are identically $N(0,\sigma_s^2)$ distributed and the $\varepsilon_{ijk}$
are identically $N(0,\sigma_{\varepsilon}^2)$ distributed,
where $\sigma_s^2 > 0$ and $\sigma_{\varepsilon}^2 > 0$.
Let $m = (1/n_1)+(1/n_2)$.

The parameter of interest is $\theta = \phi_1 - \phi_2$. The parameter describing the differential
carryover effect is $\psi= 3(\lambda_1 - \lambda_2)/4$.
Let $\bar Y_{i \boldsymbol{\cdot} k} = (1/n_i) \sum_{j=1}^{n_i} Y_{ijk}$
($i=1,2$; $k=1,2,3,4$). We reduce that data to $D_1$, $D_2$, $D_3$, $D_4$, where
$D_1 = \bar{Y}_{1 \boldsymbol{\cdot} 1} - \bar{Y}_{2 \boldsymbol{\cdot} 1}$,
$D_2 = \bar{Y}_{1 \boldsymbol{\cdot} 2} - \bar{Y}_{2 \boldsymbol{\cdot} 2}$,
$D_3 = \bar{Y}_{1 \boldsymbol{\cdot} 3} - \bar{Y}_{2 \boldsymbol{\cdot} 3}$ and
$D_4 = \bar{Y}_{1 \boldsymbol{\cdot} 4} - \bar{Y}_{2 \boldsymbol{\cdot} 4}$.
The motivation for this data reduction is presented in Appendix A. Let
\begin{equation*}
A = \frac{1}{4} (D_1 - D_2 + D_3 - D_4).
\end{equation*}
This is the usual estimator of $\theta$, when it is assumed that $\psi = 0$
(see e.g. Table I of Senn \& Lambrou, 1998).
Let
\begin{equation*}
\hat \Theta = D_1 - \frac{1}{4} D_2 - \frac{1}{2} D_3 - \frac{1}{4} D_4.
\end{equation*}
This is an unbiased estimator
(which is unaffected by differential carryover)
of $\theta$ that is unaffected by between-subject variation
(cf Table I of Senn \& Lambrou, 1998). We will use
$\hat \Theta$ as the estimator of $\theta$ when
it cannot be assumed that necessarily
$\psi=0$.
We view
$\hat \Theta$ as the analogue for an ABAB/BABA design of the estimator of $\theta$ constructed
using only data from the first period of an AB/BA design.
As shown in Appendix A, $\hat \Theta$ is inefficient by comparison
with the usual estimator of $\theta$ based on data from a completely randomized trial,
using the same number of measurements of response, unless $\sigma_s^2 \ge 4.5 \, \sigma_{\varepsilon}^2$.
We will also make use of the following unbiased estimator of $\psi$:
\begin{equation*}
\hat \Psi = \frac{3}{4} (D_1 - D_3).
\end{equation*}
As shown in Appendix B, these statistics have the following distributions:
$A \sim N(\theta - \psi, m \sigma_{\varepsilon}^2/4)$,
$\hat \Theta \sim N(\theta, 11 m \sigma_{\varepsilon}^2/8)$ and
$\hat \Psi \sim N(\psi, 9 m \sigma_{\varepsilon}^2/8)$. Note that when $\psi = 0$, $A$ is a much more
efficient estimator of $\theta$ than $\hat \Theta$.

To deal with the possibility of non-zero differential carryover, it is tempting to consider the use
of the following two-stage procedure. In the first stage a preliminary test of the null
hypothesis that $\psi=0$ (against the alternative that $\psi \ne 0$) is carried out,
using a test statistic based
on $\hat \Psi$. If this null hypothesis is accepted
then the confidence interval for $\theta$ is constructed using the estimator $A$
and having nominal coverage
%probability
$1-\alpha$,
assuming that $\psi=0$.
If, on the other hand, this null hypothesis is rejected then this confidence interval is constructed
to have nominal coverage $1-\alpha$, using
the estimator $\hat \Theta$.

To analyse the properties of this two-stage procedure, we make the
simplification that $\sigma_{\varepsilon}^2$ is known.
Freeman (1989) makes the same simplification. So, in the first stage,
we test the null hypothesis $H_0: \psi = 0$ against the alternative hypothesis $H_1: \psi \ne 0$
using the test statistic $\sqrt{8/9m} \; \hat \Psi / \sigma_{\varepsilon}$. This test statistic
has an $N(0,1)$ distribution under $H_0$. Define the quantile
$c_a$ by the requirement that $P(-c_a \le Z \le c_a) = 1 - a$
for $Z \sim N(0,1)$. The following is a test of $H_0$ against $H_1$, with level of significance
$\alpha_1$. Accept $H_0$ if $\big |\sqrt{8/9m} \; \hat \Psi / \sigma_{\varepsilon} \big| < c_{\alpha_1}$;
otherwise reject $H_0$. In the second stage we proceed as follows. If $H_0$ is accepted then we construct
a confidence interval for $\theta$, with nomimal coverage $1-\alpha$, assuming that $\psi = 0$. This confidence
interval is
\begin{equation}
\label{CI_A}
\Big [ A - c_{\alpha} \sqrt{m/4} \, \sigma_{\varepsilon}, \, A + c_{\alpha} \sqrt{m/4} \, \sigma_{\varepsilon} \Big ].
\end{equation}
If, on the other hand, $H_0$ is rejected then we do not assume that $\psi = 0$ and we construct
a confidence interval for $\theta$, with nomimal coverage $1-\alpha$, based on $\hat \Theta$. This confidence
interval is
\begin{equation}
\label{CI_Thetahat}
\Big [ \hat \Theta - c_{\alpha} \sqrt{11 m / 8} \, \sigma_{\varepsilon},
\, \hat \Theta + c_{\alpha} \sqrt{11 m / 8} \, \sigma_{\varepsilon} \Big ].
\end{equation}

Let $J$ denote the confidence interval for $\theta$ that results from this two-stage procedure.
Also let $\gamma = \sqrt{8/9m} \, \psi/\sigma_{\varepsilon}$.
As shown in Appendix A, the coverage probability of the
confidence interval $J$ is
\begin{equation}
\label{formula_cov}
P(\theta \in J) = P \big ( |H| < c_{\alpha_1} \big ) P \big ( |X| \le c_{\alpha} \big )
+ P \big ( |G| \le c_{\alpha}, |H| \ge c_{\alpha_1} \big ),
\end{equation}
where
\begin{equation}
\label{dist_G_H}
\begin{bmatrix}G \\ H\end{bmatrix}  \sim N\Bigg( \begin{bmatrix}0 \\ \gamma  \end{bmatrix} , \begin{bmatrix} 1 & 3/\sqrt{11} \, \\
 3/\sqrt{11} & 1 \end{bmatrix}\Bigg)
\end{equation}
\\
and $X \sim N(-3 \gamma / \sqrt{2}\, , \, 1)$. Note that, for given $\alpha_1$ and $\alpha$, the coverage
probability \eqref{formula_cov} is a function of the scaled differential carryover $\gamma$.
The right-hand side of \eqref{formula_cov}
is easily computed (using e.g. R or MATLAB programs),
for each given $\gamma$. The last term on the right-hand side of \eqref{formula_cov}
can be computed by evaluating the cumulative distribution function of the bivariate normal distribution
\eqref{dist_G_H}. Alternatively, this term can be computed by numerically evaluating the integral \eqref{final}, derived in
Appendix C.

\bigskip
%\newpage

\begin{center}
{\large\textbf{3. Numerical evaluation of the coverage probability as a function of $\boldsymbol{\gamma}$}}
\end{center}

Consider the two-stage procedure, for an AB/BA trial, based on a preliminary test with given level
of significance and resulting in a confidence interval with a given nominal coverage. As shown by Freeman (1989),
the actual coverage probability of this confidence interval
depends on both the scaled differential carryover $\big(\lambda \sqrt{n}/\sigma$ in Freeman's
notation$\big)$ and
$\rho = \sigma_s^2 / (\sigma_{\varepsilon}^2 + \sigma_s^2)$.
For each different value of $\rho$, there is a different graph of this coverage probability as a function of the
scaled differential carryover. The larger the value of $\rho$,
the smaller the minimum coverage probability of this confidence interval.

Now consider the two-stage procedure described in the previous section, for an ABAB/BABA trial,
based on a preliminary test with given level
of significance and resulting in a confidence interval with a given nominal coverage. In sharp contrast
to the AB/BA trial, the actual coverage probability (given by \eqref{formula_cov}) of this confidence interval depends {\sl only} on the
scaled differential carryover $\gamma$. This coverage is {\sl uninfluenced} by the between-subject variability
(which is described by the parameter $\sigma_s^2$).
For level of significance $\alpha_1 = 0.1$ of the preliminary test and nominal coverage $1-\alpha = 0.95$, this coverage
probability as a function of $\gamma$ is shown in Figure 1.
The minimum coverage probability of this confidence interval
is 0.4711, showing that this confidence interval is completely inadequate.
The minimum coverage
probability of this confidence interval was computed for a wide range of values of $\alpha_1$ and $1-\alpha$. In every case,
this confidence interval was found to have minimum coverage probability far below nominal, showing that
it is completely inadequate.
Note that for a given
level of significance $\alpha_1$ of the preliminary test and given nominal coverage $1-\alpha$, the minimum coverage probability
of this confidence interval does not depend on either of the sample sizes $n_1$ and $n_2$. The only effect of an
increase in $n_1$ and $n_2$ is to change the scaling (via $m=(1/n_1)+(1/n_2)$) of the differential carryover $\psi$.
Consequently, the harmful effect of preliminary hypothesis testing does not disappear with an increase in sample sizes
$n_1$ and $n_2$.

\newpage

%\medskip

% The following command was used to prevent the figures moving around too much.
\FloatBarrier

\begin{figure}[h]
\label{Figure1}
%\vskip-3cm
    %\includegraphics[keepaspectratio=true, width=15cm]{fig1_vilnius_13jun06}
    %\centering
  \hskip-0.7cm
    \includegraphics[scale=0.65]{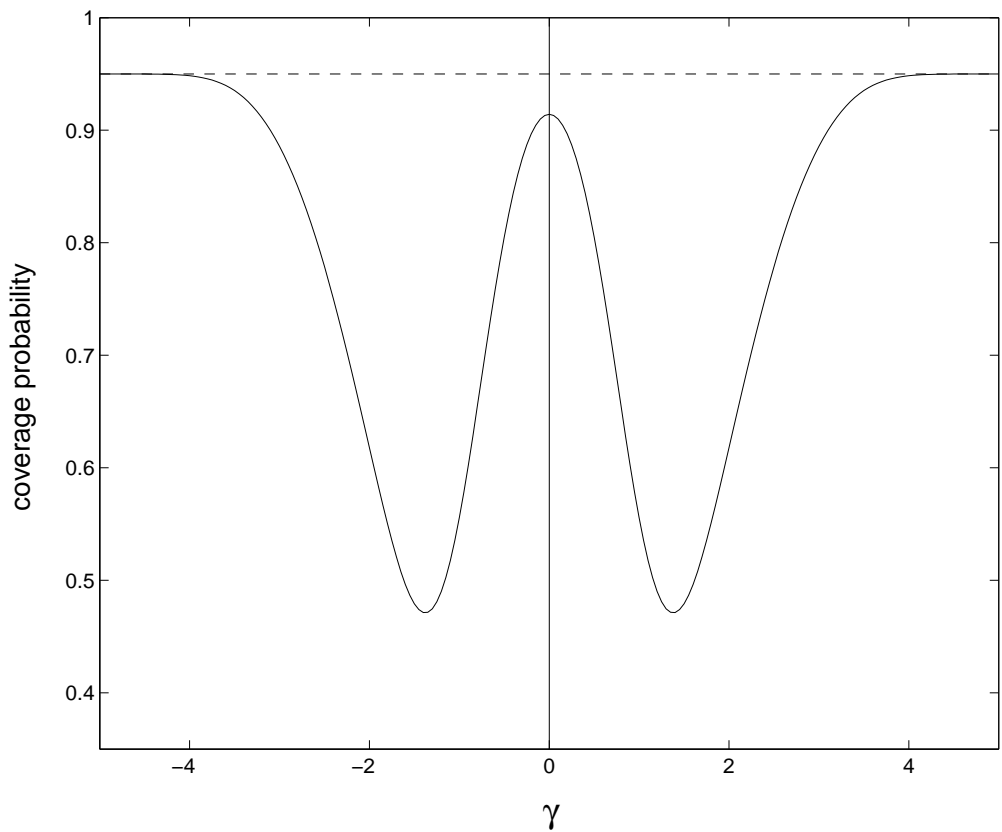}
    %\vskip-2.7cm
    \caption{Plot of the coverage
probability of the confidence interval for $\theta$, resulting from the two-stage
procedure, against $\gamma$. This confidence interval has nominal coverage $1-\alpha = 0.95$.
The preliminary hypothesis test has significance level $\alpha_1=0.1$. The horizontal dashed line has
vertical axis intercept 0.95.}
\end{figure}
\FloatBarrier

\phantom{1}
%\bigskip
\vspace{-3mm}

\begin{center}
{\large\textbf{4. Conclusion}}
\end{center}

For an ABAB/BABA trial, we have shown that the minimum coverage probability of the confidence interval resulting from
the two-stage procedure is far below the nominal coverage, showing that this confidence interval is completely inadequate.
Increasing the sample sizes $n_1$ and $n_2$ does not improve the situation.
Our conclusion is that this confidence interval should not be used.
This is similar to the conclusion of Freeman (1989) for
confidence intervals resulting from
a two-stage procedure applied to an AB/BA trial.
In other words, we provide further support for the rejection by Senn (2002, p.12) of analyses of data from any
two-treatment crossover trial
based on a preliminary test of the null hypothesis that the differential carryover is zero.

\bigskip

\begin{center}
{\large\textbf{Appendix A: The efficiency of $\boldsymbol{\hat{\Theta}}$ by comparison with
an estimator
%of $\boldsymbol{\theta}$
%based on data
from a completely randomized trial}}
\end{center}

In this appendix, we consider the efficiency of $\hat{\Theta}$ by comparison with the usual estimator of $\theta$
based on data from a completely randomized trial. For an ABAB/BABA crossover trial, the total number of
measurements of the response is $4(n_1+n_2)$. We therefore compare $\hat{\Theta}$ with the usual estimator
of $\theta$
based on data from a completely randomized trial, with $2(n_1+n_2)$ randomly-chosen subjects in each group.

Let $Y_1^A, \ldots, Y_{2(n_1+n_2)}^A$ denote the responses of the $2(n_1+n_2)$ subjects given treatment A.
Also let $Y_1^B, \ldots, Y_{2(n_1+n_2)}^B$ denote the responses of the $2(n_1+n_2)$ subjects given treatment B.
Consistently with the model \eqref{model}, we suppose that
$Y_1^A, \ldots, Y_{2(n_1+n_2)}^A,Y_1^B, \ldots, Y_{2(n_1+n_2)}^B$ are independent random variables,
where \newline
$Y_1^A, \ldots, Y_{2(n_1+n_2)}^A$ are identically $N \big(\mu+\phi_1, \sigma_{\varepsilon}^2 + \sigma_s^2 \big)$
distributed and $Y_1^B, \ldots, Y_{2(n_1+n_2)}^B$ are identically
$N \big(\mu+\phi_2, \sigma_{\varepsilon}^2 + \sigma_s^2 \big)$
distributed.
The usual estimator of $\theta$ is
\begin{equation*}
\tilde{\Theta} = \frac{1}{2(n_1+n_2)}  \big(Y_1^A + \cdots Y_{2(n_1+n_2)}^A \big) -
\frac{1}{2(n_1+n_2)}  \big(Y_1^B + \cdots Y_{2(n_1+n_2)}^A \big).
\end{equation*}
This estimator has an $N \big(\theta, (\sigma_{\varepsilon}^2 + \sigma_s^2)/(n_1+n_2) \big)$ distribution.
Suppose, for simplicity, that $n_1 = n_2 = n$. Thus
$\text{Var}(\tilde{\Theta}) = (\sigma_{\varepsilon}^2 + \sigma_s^2)/(2n)$ and
$\text{Var}(\hat{\Theta}) = 11 \sigma_{\varepsilon}^2/(4n)$. Thus
$\text{Var}(\hat{\Theta}) \le \text{Var}(\tilde{\Theta})$ if and only if
$\sigma_s^2 \ge 4.5 \, \sigma_{\varepsilon}^2$.

\bigskip

\begin{center}
{\large\textbf{Appendix B: Details for Section 2}}
\end{center}

This appendix consists of 3 sections. In the first section, we carry out data reduction. In
the second section, we derive the distributions of the statistics $A$, $\hat \Theta$ and
$\hat \Psi$. In the third section, we derive the formula \eqref{formula_cov} for the coverage probability of the
confidence interval $J$ resulting from the two-stage procedure.

\newpage

%\bigskip

\noindent \textbf{Data reduction}

\noindent It follows from the model \eqref{model} that
\begin{align*}
Y_{1j1} &= \mu + \xi_{1j} + \pi_1 + \phi_1 + \varepsilon_{1j1}\\
Y_{1j2} &= \mu + \xi_{1j} + \pi_2 + \phi_2 + \lambda_1 + \varepsilon_{1j2}\\
Y_{1j3} &= \mu + \xi_{1j} + \pi_3 + \phi_1 + \lambda_2 + \varepsilon_{1j3}\\
Y_{1j4} &= \mu + \xi_{1j} + \pi_4 + \phi_2 + \lambda_1 + \varepsilon_{1j4}\\
Y_{2j1} &= \mu + \xi_{2j} + \pi_1 + \phi_2 + \varepsilon_{2j1}\\
Y_{2j2} &= \mu + \xi_{2j} + \pi_2 + \phi_1 + \lambda_2 + \varepsilon_{2j2}\\
Y_{2j3} &= \mu + \xi_{2j} + \pi_3 + \phi_2 + \lambda_1 + \varepsilon_{2j3}\\
Y_{2j4} &= \mu + \xi_{2j} + \pi_4 + \phi_1 + \lambda_2 + \varepsilon_{2j4}
\end{align*}
Let $\bar Y_{i \boldsymbol{\cdot} k} = (1/n_i) \sum_{j=1}^{n_i} Y_{ijk}$ ($i=1,2$; $k=1,2,3,4$).
We first reduce the data to $\bar Y_{1 \boldsymbol{\cdot}  1}$, $\bar Y_{1 \boldsymbol{\cdot} 2}$,
$\bar Y_{1 \boldsymbol{\cdot} 3}$, $\bar Y_{1 \boldsymbol{\cdot} 4}$,
$\bar Y_{2 \boldsymbol{\cdot} 1}$, $\bar Y_{2 \boldsymbol{\cdot} 2}$, $\bar Y_{2 \boldsymbol{\cdot} 3}$
and $\bar Y_{2 \boldsymbol{\cdot} 4}$. Note that
\begin{align*}
\bar{Y}_{1 \boldsymbol{\cdot} 1} &= \mu + \bar{\xi}_{1 \boldsymbol{\cdot}} + \pi_1 + \phi_1
+ \bar{\varepsilon}_{1 \boldsymbol{\cdot} 1}\\
\bar{Y}_{1 \boldsymbol{\cdot} 2} &= \mu + \bar{\xi}_{1 \boldsymbol{\cdot}} + \pi_2 + \phi_2 + \lambda_1
+ \bar{\varepsilon}_{1 \boldsymbol{\cdot} 2}\\
\bar{Y}_{1 \boldsymbol{\cdot} 3} &= \mu + \bar{\xi}_{1 \boldsymbol{\cdot}} + \pi_3 + \phi_1 + \lambda_2
+ \bar{\varepsilon}_{1 \boldsymbol{\cdot} 3}\\
\bar{Y}_{1 \boldsymbol{\cdot} 4} &= \mu + \bar{\xi}_{1 \boldsymbol{\cdot}} + \pi_4 + \phi_2 + \lambda_1
+ \bar{\varepsilon}_{1 \boldsymbol{\cdot} 4}\\
\bar{Y}_{2 \boldsymbol{\cdot} 1} &= \mu + \bar{\xi}_{2 \boldsymbol{\cdot}} + \pi_1 + \phi_2
+ \bar{\varepsilon}_{2 \boldsymbol{\cdot} 1}\\
\bar{Y}_{2 \boldsymbol{\cdot} 2} &= \mu + \bar{\xi}_{2 \boldsymbol{\cdot}} + \pi_2 + \phi_1 + \lambda_2
+ \bar{\varepsilon}_{2 \boldsymbol{\cdot} 2}\\
\bar{Y}_{2 \boldsymbol{\cdot} 3} &= \mu + \bar{\xi}_{2 \boldsymbol{\cdot}} + \pi_3 + \phi_2 + \lambda_1
+ \bar{\varepsilon}_{2 \boldsymbol{\cdot} 3}\\
\bar{Y}_{2 \boldsymbol{\cdot} 4} &= \mu + \bar{\xi}_{2 \boldsymbol{\cdot}} + \pi_4 + \phi_1 + \lambda_2
+ \bar{\varepsilon}_{2 \boldsymbol{\cdot} 4}
\end{align*}
where $\bar{\xi}_{i \boldsymbol{\cdot}} = (1/n_i) \sum_{j=1}^{n_i} \xi_{ij}$ and
$\bar{\varepsilon}_{i \boldsymbol{\cdot} k} = (1/n_i) \sum_{j=1}^{n_i} \varepsilon_{ijk}$.
Note that $\bar{\xi}_{1 \boldsymbol{\cdot}}$, $\bar{\xi}_{2 \boldsymbol{\cdot}}$, $\bar{\varepsilon}_{1 \boldsymbol{\cdot} 1},
\ldots, \bar{\varepsilon}_{1 \boldsymbol{\cdot} 4}$, $\bar{\varepsilon}_{2 \boldsymbol{\cdot} 1}, \ldots,
\bar{\varepsilon}_{2 \boldsymbol{\cdot} 4}$ are independent,
$\bar{\xi}_{1 \boldsymbol{\cdot}} \sim N(0, \sigma_s^2/n_1)$,
$\bar{\xi}_{2 \boldsymbol{\cdot}} \sim N(0, \sigma_s^2/n_2)$,
$\bar{\varepsilon}_{1 \boldsymbol{\cdot} 1},
\ldots, \bar{\varepsilon}_{1 \boldsymbol{\cdot} 4}$ are identically $N(0, \sigma_{\varepsilon}^2/n_1)$ distributed
and
$\bar{\varepsilon}_{2 \boldsymbol{\cdot} 1},
\ldots, \bar{\varepsilon}_{2 \boldsymbol{\cdot} 4}$ are identically $N(0, \sigma_{\varepsilon}^2/n_2)$ distributed.

The only way to remove the influence of the parameters $\pi_1, \ldots, \pi_4$ on the
reduced data
$\bar Y_{1 \boldsymbol{\cdot}  1}, \bar Y_{1 \boldsymbol{\cdot} 2}, \ldots,
\bar Y_{2 \boldsymbol{\cdot} 4}$
is to perform
a further data reduction to $D_1, \ldots, D_4$, where
$D_1 = \bar{Y}_{1 \boldsymbol{\cdot} 1} - \bar{Y}_{2 \boldsymbol{\cdot} 1}$,
$D_2 = \bar{Y}_{1 \boldsymbol{\cdot} 2} - \bar{Y}_{2 \boldsymbol{\cdot} 2}$,
$D_3 = \bar{Y}_{1 \boldsymbol{\cdot} 3} - \bar{Y}_{2 \boldsymbol{\cdot} 3}$ and
$D_4 = \bar{Y}_{1 \boldsymbol{\cdot} 4} - \bar{Y}_{2 \boldsymbol{\cdot} 4}$.
Note that
\begin{align*}
D_1 &= \bar{\xi}_{1 \boldsymbol{\cdot}} - \bar{\xi}_{2 \boldsymbol{\cdot}} + \theta + \eta_1 \\
D_2 &= \bar{\xi}_{1 \boldsymbol{\cdot}} - \bar{\xi}_{2 \boldsymbol{\cdot}} - \theta + \frac{4}{3} \psi + \eta_2 \\
D_3 &= \bar{\xi}_{1 \boldsymbol{\cdot}} - \bar{\xi}_{2 \boldsymbol{\cdot}} + \theta - \frac{4}{3} \psi + \eta_3 \\
D_2 &= \bar{\xi}_{1 \boldsymbol{\cdot}} - \bar{\xi}_{2 \boldsymbol{\cdot}} - \theta + \frac{4}{3} \psi + \eta_4
\end{align*}
where $\eta_1 = \bar{\varepsilon}_{1 \boldsymbol{\cdot} 1} - \bar{\varepsilon}_{2 \boldsymbol{\cdot} 1}$,
$\eta_2 = \bar{\varepsilon}_{1 \boldsymbol{\cdot} 2} - \bar{\varepsilon}_{2 \boldsymbol{\cdot} 2}$,
$\eta_3 = \bar{\varepsilon}_{1 \boldsymbol{\cdot} 3} - \bar{\varepsilon}_{2 \boldsymbol{\cdot} 3}$ and
$\eta_4 = \bar{\varepsilon}_{1 \boldsymbol{\cdot} 4} - \bar{\varepsilon}_{2 \boldsymbol{\cdot} 4}$.
Note that $\eta_1, \ldots, \eta_4$ are independent and identically $N(0, m \sigma_{\varepsilon}^2)$ distributed.

\medskip
\noindent \textbf{Derivation of the distributions of the statistics $\boldsymbol{A}$, $\boldsymbol{\hat \Theta}$ and $\boldsymbol{\hat \Psi}$}

\noindent Note that
\begin{equation}
%\tag{A1}
\label{A_in_terms_of_etas}
A = \theta - \psi + \frac{1}{4} (\eta_1 - \eta_2 + \eta_3 - \eta_4).
\end{equation}
Thus $A \sim N(\theta - \psi, m \sigma_{\varepsilon}^2/4)$. Note that
\begin{equation}
%\tag{A2}
\label{Thetahat_in_terms_of_etas}
\hat \Theta = \theta + \eta_1 - \frac{1}{4} \eta_2 - \frac{1}{2} \eta_3 - \frac{1}{4} \eta_4
\end{equation}
and that
\begin{equation}
%\tag{A3}
\label{Psihat_in_terms_of_etas}
\hat \Psi = \psi + \frac{3}{4} (\eta_1 - \eta_3).
\end{equation}
It follows from \eqref{A_in_terms_of_etas} and \eqref{Psihat_in_terms_of_etas}
that $(A, \hat \Psi)$ has a bivariate normal distribution and that Cov$(A, \hat \Psi)=0$.
Thus $A$ and $\hat \Psi$ are independent random variables. It follows from
\eqref{Thetahat_in_terms_of_etas} and \eqref{Psihat_in_terms_of_etas} that
\begin{equation}
%\tag{A4}
\label{Thetahat_Psihat_dist}
\begin{bmatrix}\hat{\Theta} \\ \hat{\Psi}\end{bmatrix}  \sim
N\Bigg( \begin{bmatrix}\theta \\ \psi \end{bmatrix} , \frac{m \sigma^2_{\varepsilon}}{8} \begin{bmatrix} 11 & 9 \\ 9 & 9 \end{bmatrix}\Bigg).
\end{equation}

\medskip

\noindent \textbf{Derivation of the formula (4) for the coverage probability}

\noindent Define the event
\begin{equation*}
B = \left \{ \left | \sqrt{\frac{8}{9m}} \frac{\hat \Psi}{\sigma_{\varepsilon}} \right | < c_{\alpha_1} \right \}.
\end{equation*}
If this event occurs then $J$ is equal to \eqref{CI_A} and if $B^c$ occurs then $J$ is equal to
\eqref{CI_Thetahat}. By the law of total probability, the coverage probability $P(\theta \in J)$
is equal to
\begin{align*}
&P \big(B \cap \{ \theta \in J \} \big) + P \big(B^c \cap \{ \theta \in J \} \big) \\
&= P \left ( B \cap
\left \{ \left | \frac{(A-\theta)}{\sigma_{\varepsilon}} \sqrt{\frac{4}{m}} \right | \le c_{\alpha} \right \} \right )
+ P \left ( B^c \cap
\left \{ \left | \frac{(\hat{\Theta} - \theta)}{\sigma_{\varepsilon}} \sqrt{\frac{8}{11 m}} \right | \le c_{\alpha} \right \} \right ) \\
&= P(B) \,
P \left ( \left | \frac{(A-\theta)}{\sigma_{\varepsilon}} \sqrt{\frac{4}{m}} \right | \le c_{\alpha} \right )
+ P \left ( B^c \cap
\left \{ \left | \frac{(\hat{\Theta} - \theta)}{\sigma_{\varepsilon}} \sqrt{\frac{8}{11 m}} \right | \le c_{\alpha} \right \} \right )
\end{align*}
since $A$ and $\hat \Psi$ are independent random variables. Now define
\begin{equation*}
\gamma = \sqrt{\frac{8}{9m}} \; \frac{\psi}{\sigma_{\varepsilon}},
\quad  H = \sqrt{\frac{8}{9m}} \; \frac{\hat \Psi}{\sigma_{\varepsilon}},
\quad G = \frac{(\hat \Theta - \theta)}{\sigma_{\varepsilon}} \; \sqrt{\frac{8}{11 m}} \quad \text{and}
\quad X = \frac{(A - \theta)}{\sigma_{\varepsilon}} \; \sqrt{\frac{4}{m}}.
\end{equation*}
Thus, the coverage probability $P(\theta \in J)$ is given by \eqref{formula_cov}.
Note that it follows from \eqref{Thetahat_Psihat_dist} that the distribution of $(G,H)$ is given by
\eqref{dist_G_H}.

\bigskip

\begin{center}
{\large\textbf{Appendix C: Alternative expression for $\boldsymbol{P \big ( |G| \le c_{\alpha}, |H| \ge c_{\alpha_1} \big )}$}}
\end{center}

\medskip

\noindent In this appendix, we present an alternative expression for $P \big ( |G| \le c_{\alpha}, |H| \ge c_{\alpha_1} \big )$
that may be convenient for the computation  of the coverage probability \eqref{formula_cov}.
By the law of total probability,
\begin{equation*}
P \big ( |G| \le c_{\alpha} \big )
= P \big ( |G| \le c_{\alpha}, |H| \ge c_{\alpha_1} \big ) + P \big ( |G| \le c_{\alpha}, |H| < c_{\alpha_1} \big ),
\end{equation*}
so that
\begin{align*}
P \big ( |G| \le c_{\alpha}, |H| \ge c_{\alpha_1} \big )
&= P \big ( |G| \le c_{\alpha} \big ) - P \big ( |G| \le c_{\alpha}, |H| < c_{\alpha_1} \big ) \\
&= 1- \alpha - P \big ( |G| \le c_{\alpha}, |H| < c_{\alpha_1} \big )
\end{align*}
since $G \sim N(0,1)$.

Let $f_{G,H}(g,h)$ denote the probability density function of $(G,H)$, evaluated at $(g,h)$.
Also, let $f_{H|G}(h|g)$ denote the probability density function of $H$ conditional
on $G=g$, evaluated at $h$. Let $\phi$ denote the $N(0,1)$ probability density function.
Observe that
\begin{align}
\notag
P \big ( |G| \le c_{\alpha}, |H| < c_{\alpha_1} \big )
&= \int_{-c_{\alpha}}^{c_{\alpha}} \int_{-c_{\alpha_1}}^{c_{\alpha_1}} f_{G,H}(g,h) \, dh \, dg \\
\label{intermediate}
%\tag{B1}
&= \int_{-c_{\alpha}}^{c_{\alpha}} \int_{-c_{\alpha_1}}^{c_{\alpha_1}} f_{H|G}(h|g) \, dh \ \phi(g) \, dg
\end{align}
It follows from \eqref{dist_G_H}
that the distribution of $H$ conditional on $G=g$ is $N \big(\mu(g), v \big)$, where
$\mu(g) = \gamma + (3 g/\sqrt{11})$ and $v=2/11$. Thus \eqref{intermediate}
is equal to
\begin{equation}
%\tag{B2}
\label{final}
\int_{-c_{\alpha}}^{c_{\alpha}}
\Big ( \Phi \big(c_{\alpha_1}; \mu(g), v \big) - \Phi \big(-c_{\alpha_1}; \mu(g), v \big) \Big )
\, \phi(g) \, dg
\end{equation}
where $\Phi(x; \mu, v)$ denotes the $N(\mu, v)$ cumulative distribution function,
evaluated at $x$. The integral \eqref{final} is readily evaluated using the numerical quadrature
functions available in either R or MATLAB.

\bigskip

\begin{center}
\sl{References}
\end{center}

\rf ARMITAGE, P. \& HILLS, M. (1982).
The two-period crossover trial. {\sl Statistician}
{\bf 31}, 119--131.

\smallskip

\rf FLEISS, J.L. (1986). On multiperiod crossover trials. Letter to the editor. {\sl Biometrics}
{\bf 42}, 449--450.

\smallskip

\rf FLEISS, J.L. (1989). A critique of recent research on the two-treatment crossover design.
{\sl Controlled Clinical Trials}
 {\bf 10}, 237--243.

\smallskip

\rf FREEMAN, P. (1989). The performance of the two-stage analysis of two-treatment,
two-period crossover trials. {\sl Statistics in Medicine} {\bf 8}, 1421--1432.

\smallskip

\rf GRIEVE, A.P. (1987). A note on the analysis of the two-period crossover design when the period-treatment
interaction is significant. {\sl Biometrical Journal}  {\bf 7}, 771--775.

\smallskip

\rf GRIZZLE, A.P. (1965). The two-period change-over design and its use in clinical trials. {\sl Biometrics} {\bf 21}, 467--480.

\smallskip

\rf GRIZZLE, A.P. (1974). Corrigenda to Grizzle (1965). {\sl Biometrics} {\bf 42}, 459.

\smallskip

\rf HILLS, M. \& ARMITAGE, P. (1979). The two-period crossover clinical trial.
{\sl British Journal of Clinical Pharmacology} {\bf 8}, 7--20.

\smallskip

\rf JONES, B. \& KENWARD, M.G. (2003). {\sl Design and Analysis of Cross-Over Trials. Second Edition}. Boca Raton: Chapman
\& Hall.

\smallskip

\rf KABAILA, P. (2009). The coverage properties of confidence regions after model
selection. {\sl International Statistical Review} {\bf 77}, 405--414.

\smallskip

\rf SENN, S. \& LAMBROU, D. (1998). Robust and realistic approaches to carry-over.
{\sl Statistics in Medicine} {\bf 17}, 2849--2864.

\smallskip

\rf SENN, S. (2001). Cross-over trials in drug development: theory and practice.
{\sl Journal of Statistical Planning and Inference} {\bf 96}, 29--40.

\smallskip

\rf SENN, S. (2002). {\sl Cross-Over Trials in Clinical Research, Second Edition}. Chichester: Wiley.

\smallskip

\rf SENN, S. (2005). Misunderstandings regarding clinical cross-over trials. Letter to the editor.
{\sl Statistics in Medicine} {\bf 24}, 3675--3678.

\smallskip

\rf SENN, S. (2006). Crossover trials in Statistics in Medicine: the first `25'
years. {\sl Statistics in Medicine} {\bf 25}, 3430--3442.

\end{document}